\documentclass[11pt]{article}
\usepackage{bbm}
\usepackage{epsfig}
\usepackage{amssymb,amsmath,amscd, mdwmath}
\usepackage{graphicx,cite,url}

\begin{document}
\title{\ Two unconditional stable schemes for simulation of heat equation  on manifold using  DEC}

\author{
Zheng Xie$^1$\thanks{E-mail: lenozhengxie@yahoo.com.cn}~~~~Yujie
Ma$^2$\thanks{ E-mail: yjma@mmrc.iss.ac.cn This work is partially
supported by CPSFFP (No. 20090460102) NKBRPC (No. 2004CB318000) and
NNSFC (No. 10871170) }
\\{\small
$1.$ Center of Mathematical Sciences Zhejiang University
(310027),China}
\\ {\small $2.$ Key Laboratory of Mathematics Mechanization,}
\\ {\small  Chinese Academy of Sciences,  (100090), China}}

\date{}

\maketitle

\begin{abstract}

  To predict the heat diffusion in a given region over time,  it is
often necessary to find the numerical solution for heat equation.
With the techniques of discrete differential calculus, we propose
two unconditional stable numerical schemes  for simulation  heat
equation on space manifold and time. The analysis of their stability
and error is accomplished by the use of maximum principle.
\end{abstract}

\vskip 0.2cm \noindent {\bf Keywords: }Discrete exterior calculus,
Discrete manifold, Heat equation.

\vskip 0.2cm \noindent {\bf PACS(2010):} 44.05.+e, 02.30.Jr,
02.30.Mv, 02.40.Vh
\section{Introduction}

To investigate the predictions of heat equation it is often
necessary to approximate its solution
numerically\cite{larsson,morton}. Discrete exterior calculus (DEC)
constitutes a discrete realization of the exterior differential
forms structure, and therefore, the right framework in which to
develop a   discretization of exterior differential
system\cite{whitney,arnold,bossavit1,desbrun,leok,hiptmair,hyman,meyer}.
The operators in this system such as Hodge star, exterior
derivative, and Laplace operator can also be naturally discretized
by   DEC. The techniques of DEC can be used to construct a
conditional stable scheme for heat equation\cite{xie-ye-ma2}. The
stable condition is a very severe restriction, and implies that very
many time steps will be necessary to follow the solution over a
reasonably large time interval.

 In this
paper we shall propose two unconditional stable schemes for the
numerical solution of heat equation in space manifold and the time,
namely implicit and semi-implicit DEC schemes. The analysis of their
error is accomplished by the use of maximum principle. The methods
proposed here can be extended to problems with
 general boundary condition as explicit scheme, then to general linear and
nonlinear parabolic equations.

\section{Implicit DEC scheme for heat equation}

  For a function $\psi(x,y,z,t)$ of three
spatial variables $(x,y,z)$ and the time variable $t$, the heat
equation is
$$\rho c\dfrac{\partial \psi}{\partial t}-k\left(\dfrac{\partial^2\psi}{\partial x^2}
+\dfrac{\partial^2\psi}{\partial
y^2}+\dfrac{\partial^2\psi}{\partial z^2}\right)-\rho
Q=0,\eqno{(1)}$$ with boundary conditions where $\psi$ is the
temperature, $\rho$ is the material density,
  $c$ is the material specific
heat,
   $k$ is the thermal conductivity,
  $Q$ is the
internal heat source density.

The DEC method can approximate Laplace operator as
  $$\Delta\approx  \ast^{-1}d^{T}\ast + d^{T}\ast d,$$
where $d$ is the matrix of discrete exterior derivative, $\ast$ is
the matrix of discrete Hodge star operator. The explicit DEC scheme
for heat equation uses the forward time difference for the temporal
derivative of Eq.(1), which is  conditional stable\cite{xie-ye-ma2}.
If we need to refine the space mesh to improve the accuracy of the
solution the amount of work involved increases very rapidly, since
we shall also reduce the length of time step. Now, we show how the
use of a backward time difference gives a difference scheme which
avoids this restriction.

Given   $2$D or $3$D space manifold  with a smooth boundary, a
simplicial mesh is a tessellation of it by triangular or tetrahedra,
satisfies the condition that any two of them may intersect along a
common face, edge or vertice. For some situations, a source having
azimuthal symmetry about its axis is considered. In this case, we
only need to consider 2D triangular discrete manifold as the space.
The 3D case can also be done in the same way.

 At first, we define
some values on mesh. Take Fig.1  as an example for  2D mesh, in
which $0$, $A$,..., $F$ are vertices,    $1$, $2$,...,$6$ are the
circumcenters of triangles,    $a$, $b$,...,$f$ are the
circumcenters of edges. Denote $l_{ij}$ as the length of line
segment $(i,j)$ and  $A_{ijk}$ as the area of triangle $(i,j,k)$.
 $$
\begin{minipage}{0.99\textwidth}
\begin{center} \includegraphics[scale=0.25]{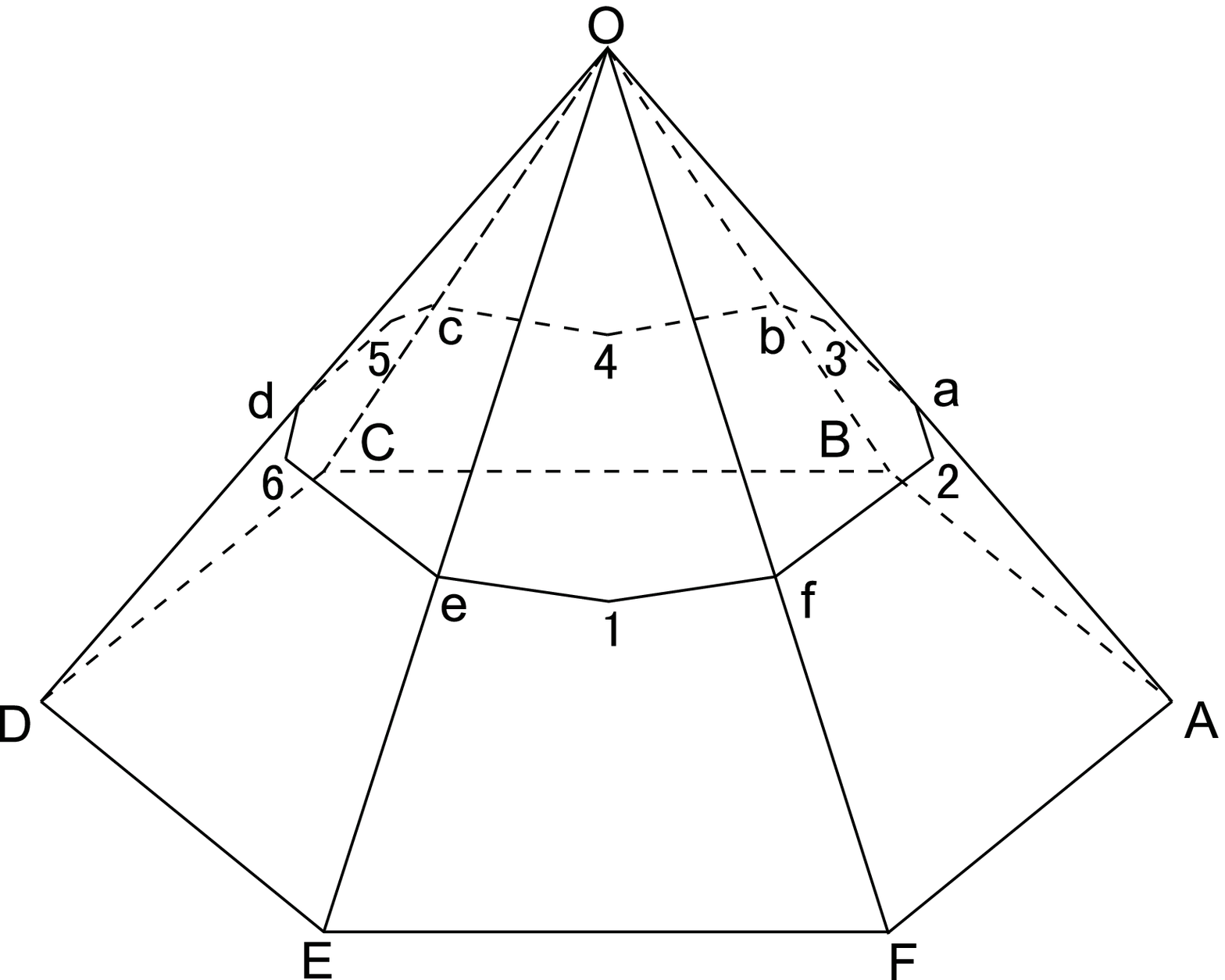}
\end{center}
\centering{ Figure.1}
\end{minipage}
$$Define $$l_{12}:=l_{1f}+l_{2f},~l_{23}:=l_{2a}+l_{3a},...,
   l_{61}:=l_{6e}+l_{1e},$$  and $$ P_{123456}:=
A_{01f}+A_{02f}+\cdots+A_{01e} .$$ The Laplace operator acting on
$\psi$ at vertice $0$   approximated by DEC is
$$ \begin{array}{lll}\Delta \psi_0&\approx& \dfrac{1}{P_{123456}} \left(\dfrac{l_{23}}{l_{A0}}
(\psi_A-\psi_0)+\dfrac{l_{34}}{l_{B0}}(\psi_B-\psi_0)+\dfrac{l_{45}}{l_{C0}}(\psi_C-\psi_0)\right.\\
&&\left.+\dfrac{l_{56}}{l_{D0}}(\psi_D-\psi_0)+\dfrac{l_{16}}{l_{E0}}(\psi_E-\psi_0)+\dfrac{l_{12}}{l_{F0}}(\psi_F-\psi_0)
\right).\end{array}
$$
If we use the backward time difference, we obtain the  implicit DEC
scheme  for Eq.(1) as follows:
$$\begin{array}{lll}  \psi^{n+1}_0 &=&\psi^n_0+\dfrac{k\Delta t}{c\rho P_{123456} }
\left(\dfrac{l_{23}
}{l_{A0}}(\psi^{n+1}_A-\psi^{n+1}_0)+\dfrac{l_{34} }{l_{B0}
}(\psi^{n+1}_B-\psi^{n+1}_0)
\right.\\
&&\left.+\dfrac{l_{45}}{l_{C0}}(\psi^{n+1}_C-\psi^{n+1}_0)
+\dfrac{l_{56}}{l_{D0}}(\psi^{n+1}_D-\psi^{n+1}_0)
+\dfrac{l_{16}}{l_{E0}}(\psi^{n+1}_E-\psi^{n+1}_0)\right.\\
&&\left.+\dfrac{l_{12}}{l_{F0}}(\psi^{n+1}_F-\psi^{n+1}_0)
\right)+\dfrac{\Delta t}{c}\rho Q^n_0.\end{array}\eqno{(2)}$$

Now, we use maximum principle to prove the unconditional stability
for scheme (2).   Given a perturbation $\varepsilon^n_i$ on each
$\psi^n_i$, the relation between $\varepsilon^n_i$
 and $\varepsilon^{n+1}_i$ can be induced from scheme (2) as
 follows:
   $$
\begin{array}{lll}
     \varepsilon^{n+1}_0
&=&\varepsilon^{n }_0-\dfrac{\Delta
t}{P_{123456}}\left(\dfrac{l_{23} }{l_{A0}}+\dfrac{l_{34} }{l_{B0}
}+\dfrac{l_{45} }{l_{C0} }+\dfrac{l_{56} }{l_{D0} }+\dfrac{l_{16}
}{l_{E0} }+\dfrac{l_{12}
}{l_{F0} }\right)  \varepsilon^{n+1 }_0 \\
&&+ \dfrac{\Delta t}{P_{123456}}\left(\dfrac{l_{23} }{l_{A0}}
\varepsilon^{n+1}_A +\dfrac{l_{34} }{l_{B0} } \varepsilon^{n +1}_B
 +\dfrac{l_{45} }{l_{C0} } \varepsilon^{n+1 }_C +\dfrac{l_{56} }{l_{D0} } \varepsilon^{n +1}_D
\right.\\
&&\left.  +\dfrac{l_{16} }{l_{E0} } \varepsilon^{n
+1}_E+\dfrac{l_{12} }{l_{F0} } \varepsilon^{n +1}_F \right).
\end{array}
$$
If we introduce the maximum perturbation at $n$ time step by writing
$$|\varepsilon^n|=\mathrm{Max}_{i\in V} |\varepsilon^n_i|,$$
 we can obtain
$$\begin{array}{lll}|\varepsilon^{n+1}_0|&\leq  \dfrac{ |\varepsilon^n|+ {
\dfrac{\Delta t}{P_{123456}}\left(\dfrac{l_{23}
}{l_{A0}}+\dfrac{l_{34} }{l_{B0} }+\dfrac{l_{45} }{l_{C0}
}+\dfrac{l_{56} }{l_{D0} }+\dfrac{l_{16} }{l_{E0} }+\dfrac{l_{12}
}{l_{F0} }\right)  } |\varepsilon^{n+1}|}{ { 1+\dfrac{\Delta
t}{P_{123456}}\left(\dfrac{l_{23} }{l_{A0}}+\dfrac{l_{34} }{l_{B0}
}+\dfrac{l_{45} }{l_{C0} }+\dfrac{l_{56} }{l_{D0} }+\dfrac{l_{16}
}{l_{E0} }+\dfrac{l_{12} }{l_{F0} }\right) } }& .
\end{array}\eqno{(3)}$$
Inequality (3) holds for all vertices in  grid, therefore
$$|\varepsilon^{n+1}|\leq |\varepsilon^n|.$$
That is to say scheme (2) is unconditional  stability.  Longer time
step  can be used than in  explicit DEC scheme.

By the definition of truncation error, the solution $\tilde{\psi}$
of  Eq.(1) satisfies the same relation as scheme (2) except for an
additional term $O((\Delta t)^2+l^2)$ on the right hand side. Thus
the error $X^n_i=\tilde{\psi}^n_i-\psi^n_i$ is determined from the
relation
 $$
\begin{array}{lll}
     X^{n+1}_0
&=& X^{n }_0-\dfrac{\Delta t}{P_{123456}}\left(\dfrac{l_{23}
}{l_{A0}}+\dfrac{l_{34} }{l_{B0} }+\dfrac{l_{45} }{l_{C0}
}+\dfrac{l_{56} }{l_{D0} }+\dfrac{l_{16} }{l_{E0} }+\dfrac{l_{12}
}{l_{F0} }\right) X^{n+1}_0 \\
&&+ \dfrac{\Delta t}{P_{123456}}\left(\dfrac{l_{23} }{l_{A0}}
X^{n+1}_A +\dfrac{l_{34} }{l_{B0} } X^{n+1}_B
 +\dfrac{l_{45} }{l_{C0} } X^{n+1}_C +\dfrac{l_{56} }{l_{D0} } X^{n+1}_D \right.
\\
&&\left. +\dfrac{l_{16} }{l_{E0} }
 X^{n+1}_E+\dfrac{l_{12} }{l_{F0} } X^{n+1}_F \right)+O((\Delta
t)^2+l^2)
\end{array}
$$
Define   $$|X^n|=\mathrm{Max}_{i\in V} |X^n_i|.$$ It follows that
$$\begin{array}{lll}|X^{n+1}_0|&\leq  & |X^n|+O((\Delta t)^2+l^2),
\end{array} $$ and hence that
$$|X^{n+1}|\leq |X^n|+O((\Delta t)^2+l^2),$$
so that, since $X^0=0$,
$$|X^{n+1}|\leq O((\Delta t)^2+l^2), $$
and this tends to $0$ along the refinement path under the assumed
hypotheses. So far we have assumed that numerical errors arise from
the truncation errors of the finite difference approximations, but
that the boundary values are used exactly.

\section{Semi-implicit DEC scheme for heat equation}

Implicit DEC scheme  is not so easy to use as the explicit DEC
scheme. We must solve a system of equations to give the values
simultaneously, which will cost a lot of  time. Now we propose a
semi-implicit scheme, which combines the virtues of explicit scheme
(direct calculation) and implicit scheme (unconditional stability).
This scheme is that the calculated element using next point of time
whereas the other elements  using current point of time. Hence the
semi-implicit scheme of Eq.(1)   is
$$\begin{array}{lll}  \psi^{n+1}_0 &=&\psi^n_0+\dfrac{k\Delta t}{c\rho|P_{123456}|}
\left(\dfrac{l_{23} }{l_{A0}}(\psi^n_A-\psi^{n+1}_0)+\dfrac{l_{34}
}{l_{B0} }(\psi^n_B-\psi^{n+1}_0)
\right.\\
&&\left.+\dfrac{l_{45}}{l_{C0}}(\psi^n_C-\psi^{n+1}_0)
+\dfrac{l_{56}}{l_{D0}}(\psi^n_D-\psi^{n+1}_0)
+\dfrac{l_{16}}{l_{E0}}(\psi^n_E-\psi^{n+1}_0)\right.\\
&&\left.+\dfrac{l_{12}}{l_{F0}}(\psi^n_F-\psi^{n+1}_0)
\right)+\dfrac{\Delta t}{c}\rho Q^n_0. \end{array}$$  It can be
rewritten as
 $$
\begin{array}{lll}
   \psi^{n+1}_0 &=&\frac{\psi^{n }_0
 + \frac{k\Delta t}{c\rho P_{123456}}\left(\frac{l_{23} }{l_{A0}}
\psi_A +\frac{l_{34} }{l_{B0} } \psi^{n }_B
 +\frac{l_{45} }{l_{C0} } \psi^{n }_C +\frac{l_{56} }{l_{D0} } \psi^{n }_D
 +\frac{l_{16} }{l_{E0} } \psi^{n }_E
 +\frac{l_{12} }{l_{F0} } \psi^{n }_F \right)+
 \frac{\Delta t}{c}\rho Q^n_0}{ 1+ \frac{\Delta t}{P_{123456}}\left(\frac{l_{23}
}{l_{A0}}+\frac{l_{34} }{l_{B0} }+\frac{l_{45} }{l_{C0}
}+\frac{l_{56} }{l_{D0} }+\frac{l_{16} }{l_{E0} }+\frac{l_{12}
}{l_{F0} }\right) }.
\end{array}\eqno{(4)}
$$

Now, we  consider the  stability and convergence for scheme (4). The
relation between perturbation $\varepsilon^n_i$
 and $\varepsilon^{n+1}_i$ can be induced from scheme (4) as
 follows:
    $$
\begin{array}{lll}
  \varepsilon^{n+1}_0 =\dfrac{\varepsilon^{n }_0 + \dfrac{k\Delta t}{c\rho P_{123456}}\left(\dfrac{l_{23} }{l_{A0}}
\varepsilon^{n }_A +\dfrac{l_{34} }{l_{B0} } \varepsilon^{n }_B
 +\dfrac{l_{45} }{l_{C0} } \varepsilon^{n }_C +\dfrac{l_{56} }{l_{D0} } \varepsilon^{n }_D
 +\dfrac{l_{16} }{l_{E0} } \varepsilon^{n }_E
 +\dfrac{l_{12} }{l_{F0} } \varepsilon^{n }_F \right)}{ 1+ \dfrac{k\Delta t}{c\rho P_{123456}}\left(\dfrac{l_{23}
}{l_{A0}}+\dfrac{l_{34} }{l_{B0} }+\dfrac{l_{45} }{l_{C0}
}+\dfrac{l_{56} }{l_{D0} }+\dfrac{l_{16} }{l_{E0} }+\dfrac{l_{12}
}{l_{F0} }\right) }.&&
\end{array}\eqno{(5)}
$$
Since the coefficients are positive, we can omit the modulus signs
in the Eq.(5)  to give inequality
$$\begin{array}{lll}|\varepsilon^{n+1}_0|&\leq  &\dfrac{ 1+\frac{k\Delta t}{c\rho P_{123456}}\left(\frac{l_{23}
}{l_{A0}}+\frac{l_{34} }{l_{B0} }+\frac{l_{45} }{l_{C0}
}+\frac{l_{56} }{l_{D0} }+\frac{l_{16} }{l_{E0} }+\frac{l_{12}
}{l_{F0} }\right)  }{ 1+\frac{k\Delta t}{c\rho
P_{123456}}\left(\frac{l_{23} }{l_{A0}}+\frac{l_{34} }{l_{B0}
}+\frac{l_{45} }{l_{C0} }+\frac{l_{56} }{l_{D0} }+\frac{l_{16}
}{l_{E0} }+\frac{l_{12} }{l_{F0} }\right) }|\varepsilon^{n }| \\
 &\leq  &| \varepsilon^n|.
\end{array}\eqno{(6)}$$
 Inequality (6) holds for all vertices in  mesh, therefore
$$|\varepsilon^{n+1}|\leq |\varepsilon^n|.$$
That is to say scheme(5) is unconditional  stability.

The   relation of error $X^n_i$ can derived form scheme (4) as
follows:
$$
\begin{array}{lll}
  X^{n+1}_0 =&\frac{X^{n }_0 + \frac{k\Delta t}{c\rho P_{123456}}\left(\frac{l_{23} }{l_{A0}}
X^n_A +\frac{l_{34} }{l_{B0} } X^{n }_B
 +\frac{l_{45} }{l_{C0} } X^{n }_C +\frac{l_{56} }{l_{D0} } X^{n }_D
 +\frac{l_{16} }{l_{E0} } X^{n }_E
 +\frac{l_{12} }{l_{F0} } X^{n }_F \right)}{ 1+ \frac{k\Delta t}{c\rho P_{123456}}\left(\frac{l_{23}
}{l_{A0}}+\frac{l_{34} }{l_{B0} }+\frac{l_{45} }{l_{C0}
}+\frac{l_{56} }{l_{D0} }+\frac{l_{16} }{l_{E0} }+\frac{l_{12}
}{l_{F0} }\right) }&
\\&+O((\Delta t)^2+l^2).&
\end{array}
$$
Because of the nonnegative coefficients, it follows that
$$|X^{n+1}_0| \leq   |X^n|+O((\Delta t)^2+l^2), $$ and hence that
$$|X^{n+1}|\leq |X^n|+O((\Delta t)^2+l^2),$$
so that, since $X^0=0$,
$$|X^{n+1}|\leq O((\Delta t)^2+l^2), $$
and this tends to $0$ along the refinement path under the assumed
hypotheses. That is to say scheme (5) is convergence.

In scheme (2) and (4), the derivative   is approximated by first
order difference. Equivalently, $\psi$ is approximated by linear
interpolation functions. Consulting the
  definition about accuracy of finite volume method, we
can also say that scheme (2) and (4) have first order temporal and
spacial accuracy.

\section{ Implementation of semi-implicit DEC scheme}
 The
implementation of semi-implicit scheme (4) for heat  equation is the
same as the explicit scheme  consisting of the following steps:
\begin{itemize}
  \item [1.]Set the simulation parameters. These
are the dimensions of the computational mesh and the size of the
time step, etc.;

 \item [2.] Initialize the mesh indexes.

  \item [3.]Assign    source.

    \item [4.]Assign boundary conditions.

  \item [5.]Compute the value of all spatial nodes  and temporarily store the
  result  in the circular buffer for further computation.

  \item [6.]Visualize the currently computed grid of spatial nodes.

   \item [7.]Repeat the  process  from the step 3, until reach the desired total number
   of iterations.
\end{itemize}
 In  the
following  examples, the parameters is $\rho=k=c=1$ and  the heat
source is $\sqrt{t/500}$.  The examples in Fig.2 and Fig.3  show
that the diffusion of heat on surfaces of dragon and torus.  The
temperature $\psi>1$ in red domain, $0<\psi\leq1$ in yellow domain,
and $\psi=0$ in blue domain.
$$
\begin{minipage}{0.99\textwidth}
\begin{center} \includegraphics[scale=0.34]{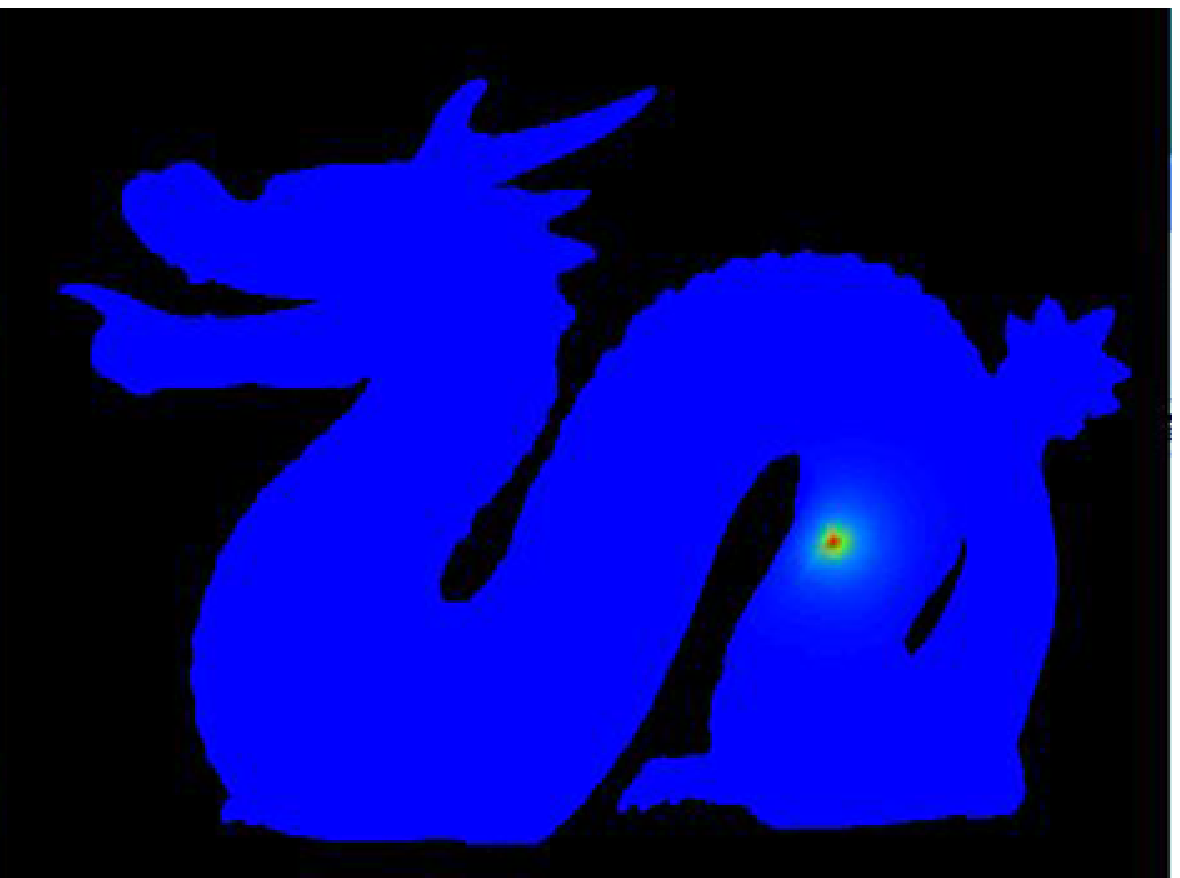}
 \includegraphics[scale=0.34]{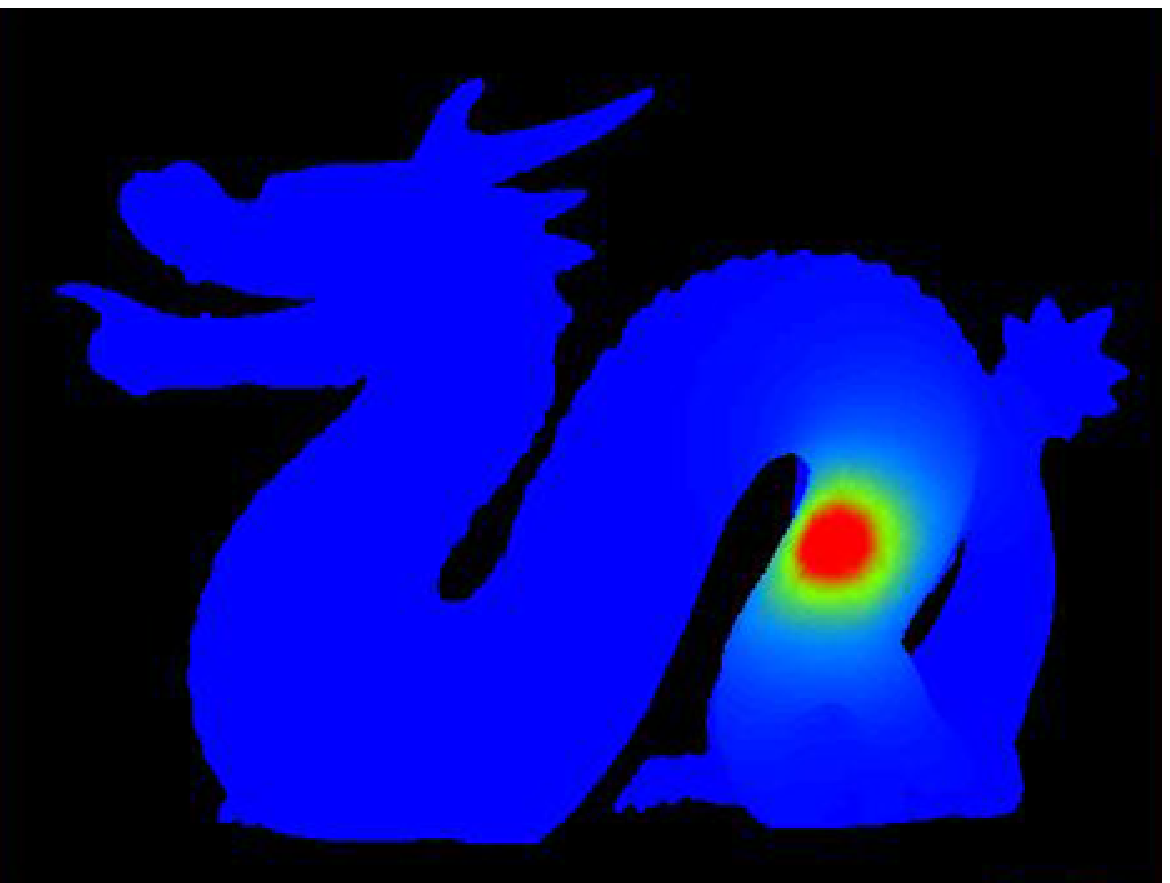}
 \includegraphics[scale=0.34]{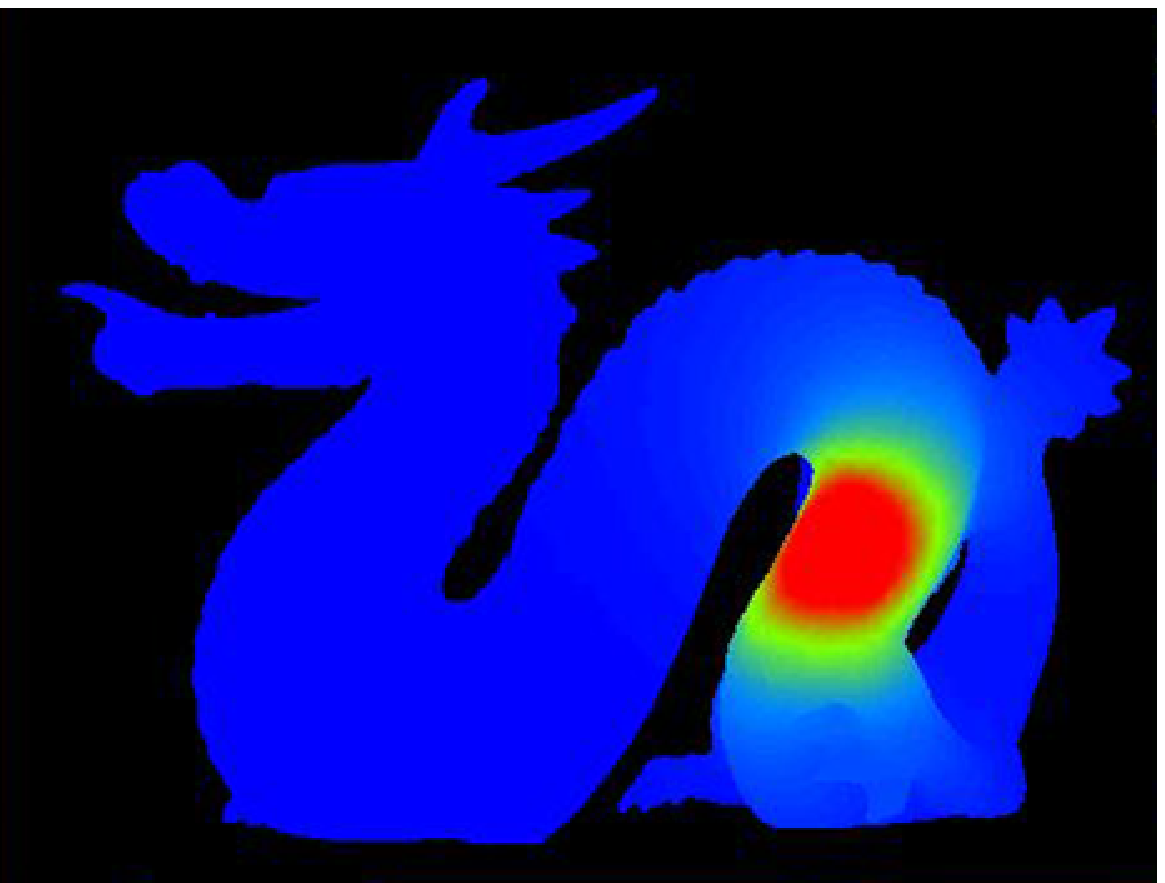}
\end{center}
\centering{ Figure.2 }
\end{minipage}
$$
$$
\begin{minipage}{0.99\textwidth}
\begin{center} \includegraphics[scale=0.35]{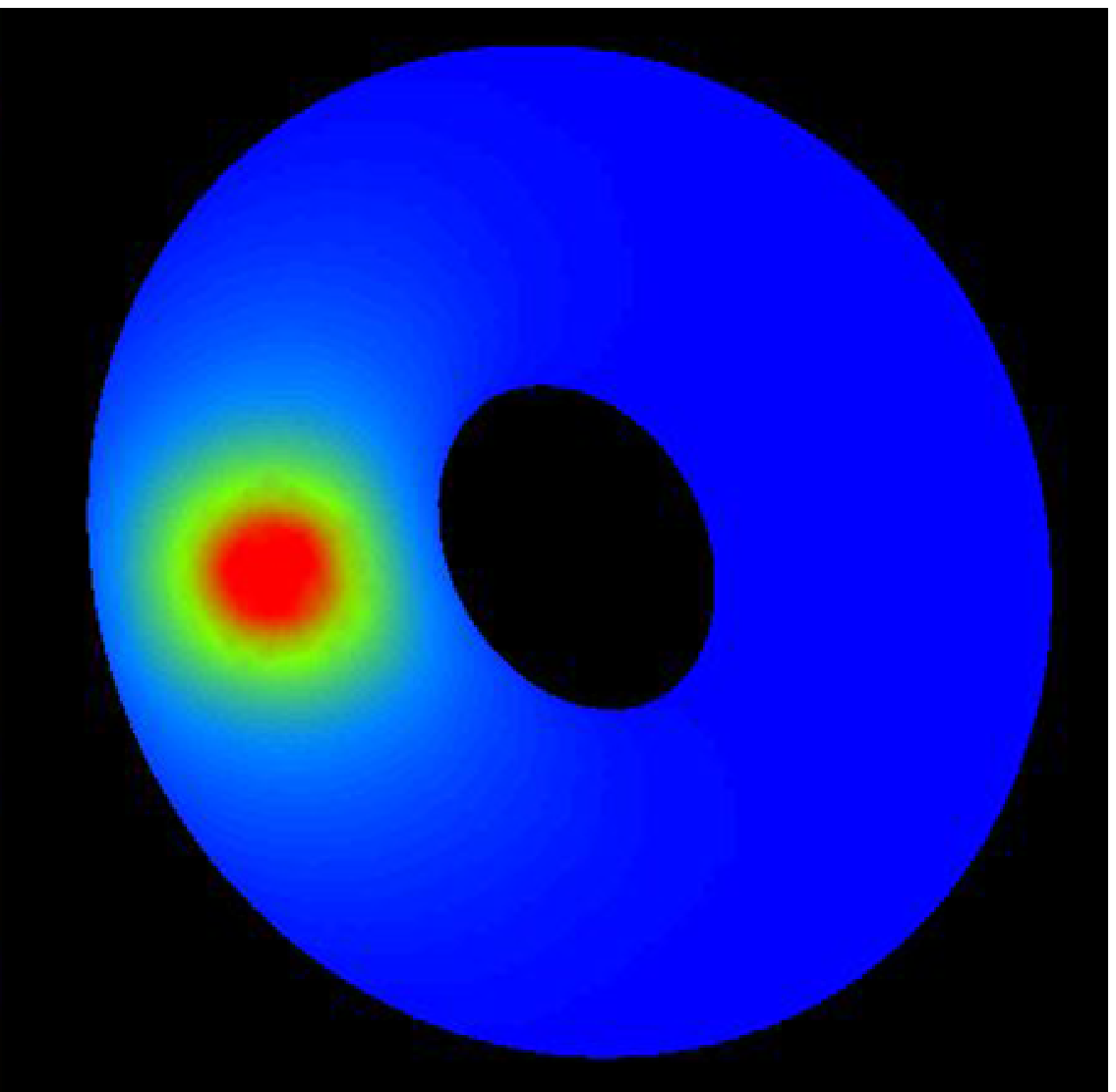}
 \includegraphics[scale=0.35]{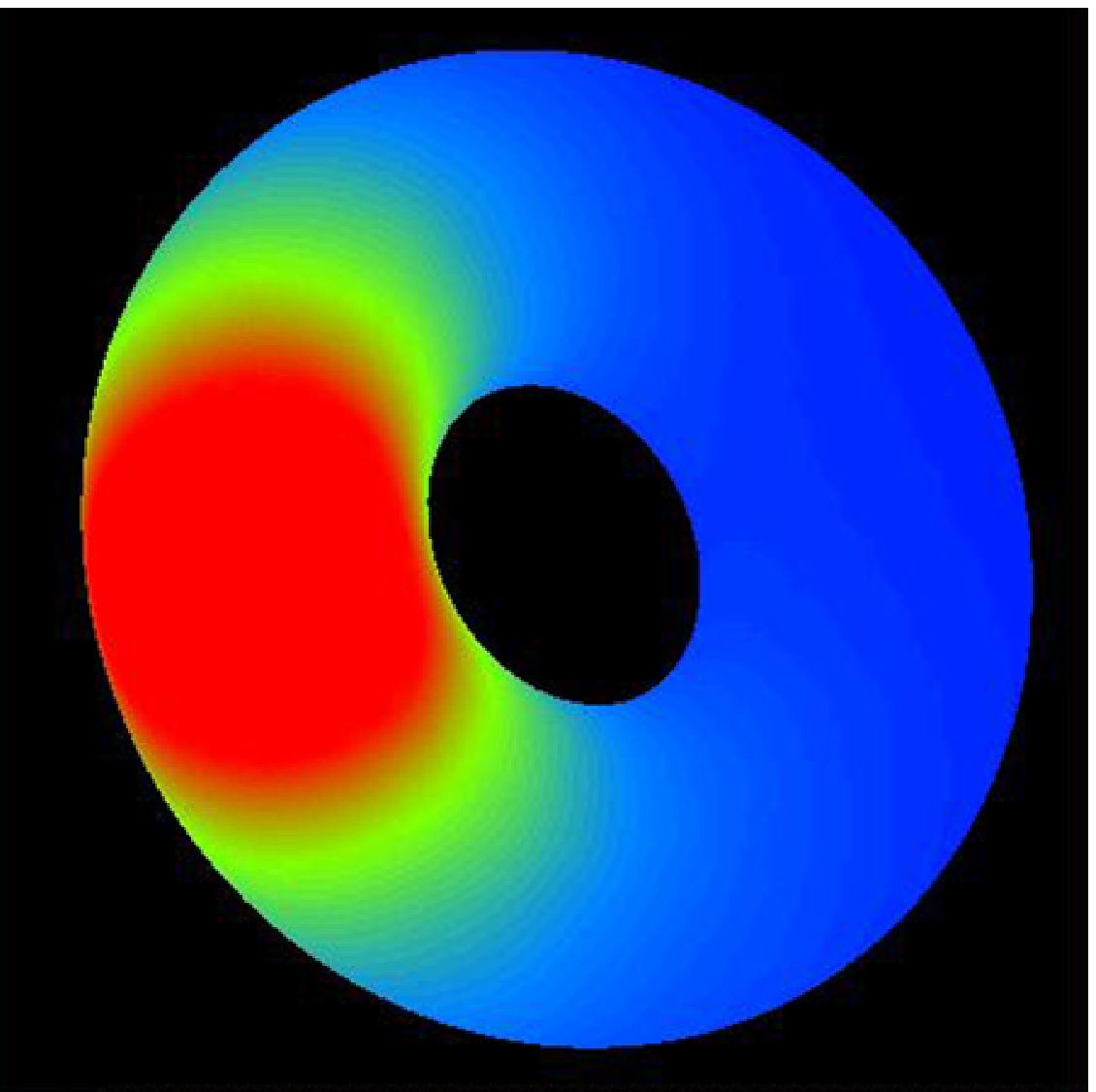}
 \includegraphics[scale=0.35]{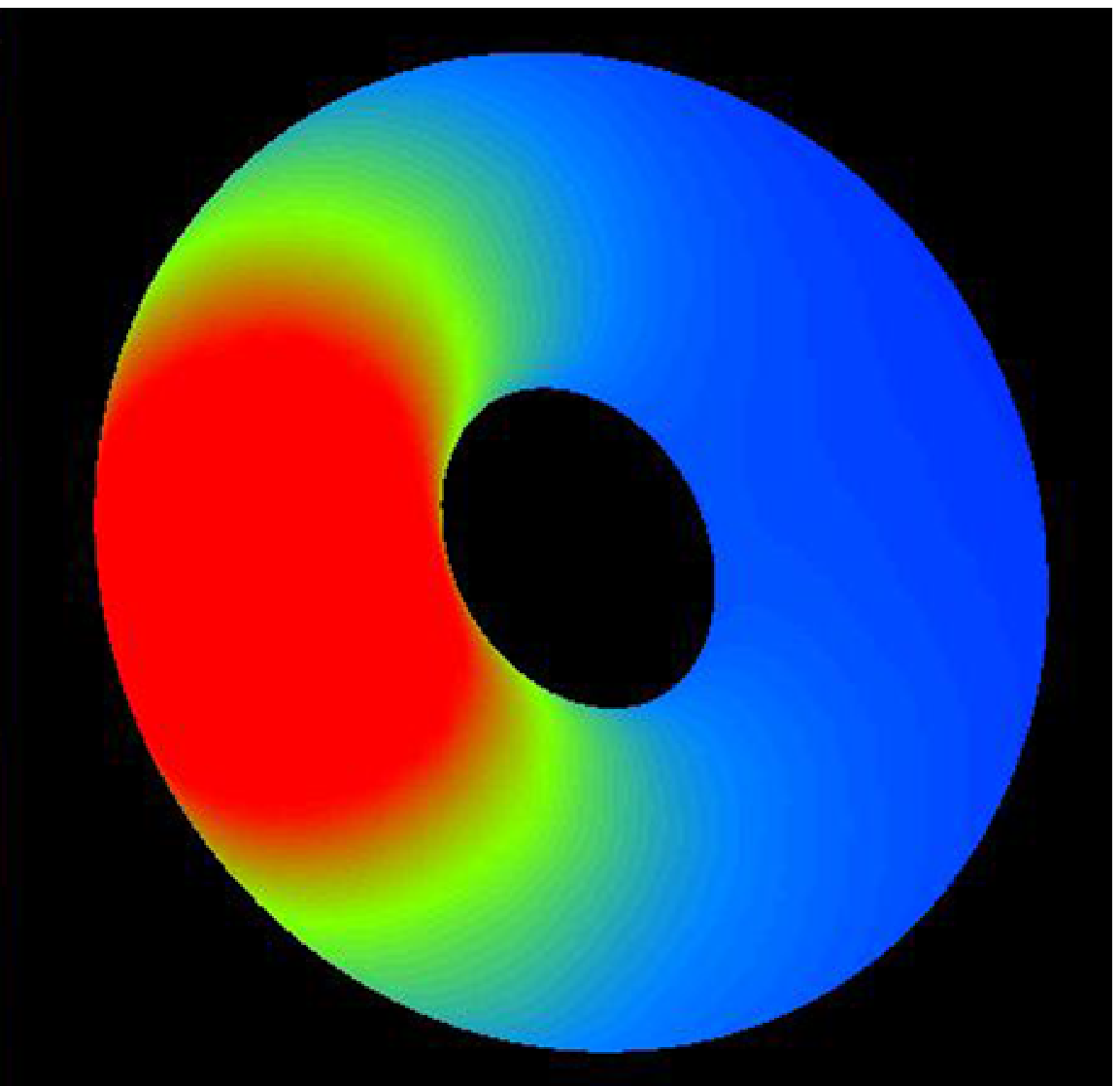}
\end{center}
\centering{ Figure.3 }
\end{minipage}
$$

\end{document}